\documentclass[12pt,twoside,reqno]{amsart}
\usepackage{amsxtra,amscd}
\usepackage{graphicx}
\usepackage{amsmath}
\usepackage{amsfonts}
\usepackage{amssymb}

\theoremstyle{definition}

\numberwithin{equation}{section}

\begin{document}
\date{2005-9-30}
\title{The Specializations in a Scheme}
\author{Feng-Wen An}
\address{School of Mathematics and Statistics, Wuhan University, Wuhan,
Hubei 430072, People's Republic of China} \email{fwan@amss.ac.cn}
\begin{abstract}
In this paper we will obtain some further properties for specializations in a
scheme. Using these results, we can take a picture for a scheme and a
picture for a morphism of schemes. In particular, we will prove that every
morphism of schemes is specialization-preserving and of norm not greater
than one (under some condition); a necessary and sufficient condition will
be given for an injective morphism between irreducible schemes.
\end{abstract}

\maketitle

\bigskip

\bigskip

\section*{Introduction}

\bigskip

Specializations are concrete and intuitive for one to study classical
varieties$^{\left[ 6\right] }$. The results on this topic relating to schemes
are mainly presented in Grothendieck's EGA. In this paper we try to obtain
some properties for the specializations in a scheme such as the lengths of
specializations.

Together with specializations, a scheme can be regarded as a partially
ordered set. In \S 1 we will prove that \emph{every morphism of schemes is
specialization-preserving} (\emph{Proposition 1.3}) and that \emph{every
specilization in a scheme is contained in an affine open subset} (\emph{%
Proposition 1.9}). Using those results, we can take a picture for a scheme (%
\emph{Remark 1.10}):

\emph{A scheme can be described to be a number of trees standing on the
ground such that}

\emph{each irreducible component is a tree;}

\emph{the generic point of an irreducible component is the root of the
corresponding tree;}

\emph{the closed points of an irreducible component are the top leaves of
the corresponding tree;}

\emph{each specialization in an irreducible component are the branches of
the corresponding tree.}

In \S 2 we will discuss the lengths of specializations, where we will notice
that the length and the dimension of a subset in a scheme are not equal in
general (Remarks 2.2-3).

Using the lengths, in \S 3 we will define the norm of a morphism of schemes
and demonstrate that \emph{any morphism of schemes is of norm not
greater than one under some condition.} Then we will obtain a picture of
morphisms of schemes (Remarks 3.6-7): \emph{As schemes are trees,
morphisms exactly scale down the trees under that condition.} A necessary
and sufficient condition will be given for an injective morphism between
irreducible schemes.

In \S 4, last section, we will present an application of specializations.

\bigskip

\textbf{Acknowledgment} The author would like to express his sincere
gratitude to Professor Li Banghe for his invaluable advice and instructions
on algebraic geometry and topology.

He also thanks Mr Yuji Odaka (Tokyo University) for pointing out the errors
in an earlier version of the preprint.

\bigskip

\bigskip

\section{Preliminaries}

\bigskip

In the section we will fix the notations and then obtain the basic facts for
specializations in a scheme. To start with, we will discuss the specializations
in a topological space (which are not Hausdorff in general) since a scheme
itself is a space.

\bigskip

Let $E$ be a topological space. Given any $x,y\in E.$ Then $y$ is a
\textbf{{{specialization}}} of $x$ (or $x$ is a \textbf{generalization} of
$y$) in $E$ if $y$ is in the closure $\overline{\{x\}}$, and we denote it by
$x\rightarrow y$ (\text{in } $E$). For $x\in E,$ we put
\[
Sp\left(  x\right)  =\{y\in E\mid x\rightarrow y\}
\]
and
\[
Gen\left(  x\right)  =\{y\in E\mid y\rightarrow x\}.
\]

If $x\rightarrow y$ and $y\rightarrow x$ both hold in $E$, $y$ is called a
\textbf{{{generic specialization}}} of $x$ in $E,$ and we denote it by
$x\leftrightarrow y$ (\text{in } $E$). The point $x$ is {{\textbf{initial}}}
if we have $x\leftrightarrow z$ for any $z\in E$ such that $z\rightarrow x;$
$x$ is {{\textbf{final}}} if we have $x\leftrightarrow z$ for any $z\in E$
such that $x\rightarrow z.$

Let $x\rightarrow y$ in $E.$ Then $y$ is said to be a {{\textbf{closest
specialization}}} of $x$ in $E$ if we have either $z=x$ or $z=y$ for any $z\in
E$ such that $x\rightarrow z$ and $z\rightarrow y$ in $E$.

\bigskip Obviously, we have the following statements:

$\left(  i\right)  $\emph{ Let }$X$\emph{ be a scheme, and }$x\in X.$\emph{
Then we have }%
\[
Sp\left(  x\right)  =\overline{\{x\}}%
\]
\emph{ and }%
\[
Gen\left(  x\right)  \cong Spec\left(  \mathcal{O}_{x}\right)  .
\]

$\left(  ii\right)  $\emph{ Let }$E$\emph{ be a topological space. Take any
}$x,y\in E.$\emph{ Then }$x\rightarrow y$\emph{ in }$E$\emph{ if and only if
}$Sp\left(  x\right)  \supseteq Sp\left(  y\right)  ;$\emph{ }%
$x\leftrightarrow y$\emph{ in }$E$\emph{ if and only if }$Sp\left(  x\right)
=Sp\left(  y\right)  $\emph{;} $x$\emph{ is final if }$x$\emph{ is a closed
point in }$E$\emph{; }$x\in E$\emph{ is initial if and only if }$x$\emph{ is a
generic point of an irreducible component of }$E$\emph{. In particular,
}$Sp\left(  x\right)  $\emph{ is an irreducible closed subset in }$E$\emph{.}

\bigskip\textbf{Example 1.1. }\emph{Let }$X$\emph{ be an Artinian scheme. Then
every point }$x\in X$\emph{ is both initial and final.}

\textbf{Proof. }As every $x\in X$ is closed, we have $Sp\left(  x\right)
=\{x\}$ for any $x\in X.\hspace{2em}\square$

\bigskip\textbf{Definition 1.2. }\emph{Let }$f:E\rightarrow F$\emph{ be a
mapping of topological spaces. }

$\left(  i\right)  $\emph{ }$f$\emph{ is said to be }$IP$\textbf{ }$-$\textbf{
preserving}\emph{ if the condition is satisfied: }

\qquad\emph{Given any closed subset }$U$\emph{ of }$E.$ \emph{Then }$f\left(
x_{0}\right)  $\emph{ is an initial point of }$\overline{f\left(  U\right)  }%
$\emph{ if }$x_{0}$\emph{ is initial in }$U.$\emph{ }

$\left(  ii\right)  $\emph{ }$f$\emph{ is \textbf{specialization-preserving}
if we have }$f\left(  x\right)  \rightarrow f\left(  y\right)  $\emph{ in }%
$F$\emph{ for any }$x\rightarrow y$\emph{ in }$E$\emph{.}

Now we obtain the main result in the section.

\bigskip{\textbf{Proposition 1.3.}} $\left(  {i}\right)  $\emph{ Every
morphism of schemes is }$IP-$\emph{preserving.}

$\left(  {ii}\right)  $\emph{ Every morphism of schemes is
specialization-preserving.}

\textbf{Proof.} It is immediate from Lemmas 1.6-7.\qquad$\square$

\bigskip

\textbf{Lemma 1.4. }\emph{Given any scheme }$X$\emph{.}

$\left(  i\right)  $\emph{ Let }$X=Spec\left(  A\right)  $\emph{ be affine.
Then we have }$x\rightarrow y$\emph{ in }$X$\emph{ for any }$x,y\in X$\emph{
if and only if }$j_{x}\subseteq j_{y}$\emph{ in }$A,$\emph{ where }$j_{x}%
$\emph{ and }$j_{y}$\emph{ denote the prime ideals in }$A$ \emph{corresponding
to }$x$\emph{\ and }$y$\emph{, respectively.}

$\left(  ii\right)  $\emph{ Take any }$x,y\in X.$\emph{ Then we have
}$x\leftrightarrow y$\emph{ in }$X$\emph{ if and only if }$x=y$\emph{.}

\textbf{Proof.} $\left(  i\right)  $ Let $x\rightarrow y$ in $X.$ We have
$Sp\left(  x\right)  =V\left(  j_{x}\right)  $ and $Sp\left(  y\right)
=V\left(  j_{y}\right)  ;$ then $Sp\left(  x\right)  \supseteq Sp\left(
y\right)  $, and hence $j_{x}\subseteq j_{y}.$ Evidently, the converse is true.

$\left(  ii\right)  $ It suffices to prove $\Rightarrow$. Let
$x\leftrightarrow y.$ Take an affine open subset $U$ of $X$ such that $x\in
U$. As $x\leftrightarrow y,$ there is the identity ${Sp}(x)={Sp}(y)$; then
\[
{Sp}(x)\bigcap U={Sp}(y)\bigcap U,
\]
which are open subsets in ${Sp}(x)$; as $Sp\left(  x\right)  $ and $Sp\left(
y\right)  $ are irreducible, we have $x,y\in U$, and hence $x\leftrightarrow
y$ in $U.$ By $\left(  i\right)  $ we have $x=y$ in $U.\qquad\square$\bigskip

Let $X_{0}$ be an irreducible closed subset of a scheme $X.$ Take an affine
open subset $U$ of $X$ such that $U\cap X_{0}\neq\varnothing$, where
$U=Spec\left(  A\right)  .$ Then
\[
U_{0}=U\cap X_{0}%
\]
is closed in $U$ and is open in $X_{0}.$ Let
\[
\Sigma=\{j_{z}\mid z\in U_{0}\}.
\]

With inclusion $\subseteq$, $\Sigma$ is a partially ordered set. There exist
minimal elements in $\Sigma$, and we denote by $\Sigma_{0}$ the set of such
minimal elements in $\Sigma.$ Let $z_{0}\in U_{0}$ with $j_{z_{0}}\in
\Sigma_{0}.$ Then $z_{0}$ is an initial point in $U_{0}.$ As $X_{0}$ is
irreducible, $U_{0}$ is irreducible; hence, $z_{0}$ is an initial point in
$X_{0}.$

If $z_{0}^{\prime}$ is another initial point in $X_{0},$ we have
$z_{0}\leftrightarrow z_{0}^{\prime}$ in $X_{0},$ and then $z_{0}%
=z_{0}^{\prime}.$ This proves (Lemma 1.6) that any irreducible closed subset
of a scheme has one and only one initial point (i.e., generic point). In
general, there is the following definition.

\bigskip

\textbf{Definition 1.5. }\emph{Assume }$E$\emph{ is a topological space
satisfying the condition: }

\qquad\emph{There exists one and only one initial point }$x_{U}$\emph{ in
every irreducible closed subset }$U$\emph{ of }$E,$\emph{ and we have }%
$x_{U}\neq x_{V}$\emph{ for any irreducible closed subset }$V$\emph{ of }%
$E$\emph{ such that }$U\neq V.$\emph{ }

\emph{Then the space }$E$\emph{ is said to have the }$\left(  UIP\right)
-$\textbf{ property}\emph{.}

\bigskip Obviously, there are many spaces which are of the $\left(
UIP\right)  -$property such as Hausdorff spaces.

\bigskip

\textbf{Lemma 1.6.} \emph{An irreducible }$T_{0}-$\emph{space which has an
initial point has the }$\left(  UIP\right)  -$\emph{property. In particular,
every scheme is of the }$(UIP)-$\emph{property.}

\bigskip

\textbf{Lemma 1.7.} \emph{Let }$f:E\rightarrow F$\emph{ be a mapping of
topological spaces.}

$\left(  i\right)  $\emph{ }$f$\emph{ is specialization-preserving if and only
if }$f$\emph{ is }$IP-$\emph{preserving.}

$\left(  ii\right)  $\emph{ Let }$F$\emph{ be of the }$\left(  UIP\right)
-$\emph{property. Then }$f$\emph{ is specialization-preserving if }$f$\emph{
is continuous.}

\textbf{Proof.} $\left(  i\right)  $ It is immediate from definition.

$\left(  ii\right)  $ Let $f$ be continuous. Take any $x\rightarrow y$ in $E.$
It is clear that $f\left(  Sp\left(  x\right)  \right)  $ is irreducible, and
then $\overline{f\left(  Sp\left(  x\right)  \right)  }$ is irreducible in
$F;$ as $f\left(  x\right)  \in f\left(  Sp\left(  x\right)  \right)  ,$ there
is
\[
Sp\left(  f\left(  x\right)  \right)  \subseteq\overline{f\left(  Sp\left(
x\right)  \right)  };
\]
as $F$ has the $\left(  UIP\right)  -$property, there is
\[
\overline{f\left(  Sp\left(  x\right)  \right)  }=Sp\left(  f\left(  x\right)
\right)  .
\]
Hence, we have
\[
\overline{f\left(  Sp\left(  y\right)  \right)  }=Sp\left(  f\left(  y\right)
\right)  .
\]

As $Sp\left(  x\right)  \supseteq Sp\left(  y\right)  ,$ we have
\[
f\left(  Sp\left(  x\right)  \right)  \supseteq f\left(  Sp\left(  y\right)
\right)  ;
\]
then there is
\[
\overline{f\left(  Sp\left(  x\right)  \right)  }\supseteq\overline{f\left(
Sp\left(  y\right)  \right)  },
\]
and it follows that
\[
Sp\left(  f\left(  x\right)  \right)  \supseteq Sp\left(  f\left(  y\right)
\right)
\]
holds. Hence, there is the specialization $f\left(  x\right)  \rightarrow
f\left(  y\right)  $ in $F.\qquad\square$

\bigskip

\textbf{Lemma 1.8.} \emph{Let }$X$\emph{ be a scheme. For every specialization
}$x\rightarrow y$\emph{ in }$X,$\emph{ there is an affine open subset }%
$U$\emph{ of }$X$\emph{ such that }$x,y\in U.$

\textbf{Proof.} Let $x\rightarrow y$ in $X.$ Hypothesize that there is no
affine open set $W$ such that $x,y\in W.$ Let $V$ be an affine open set such
that $y\in V$ but $x\not \in V.$ Then $y$ is not a limit point of the set
$\{x\},$ and we will obtain a contradiction.\qquad$\square$

\bigskip

\textbf{Proposition 1.9.} \emph{Let }$X$\emph{ be a scheme, and }$x,y\in
X$\emph{ such that }$x\rightarrow y.$\emph{ Then there is an affine open set
}$U=Spec\left(  A\right)  $\emph{ in }$X$\emph{ such that }%
\[
x,y\in U\text{ \emph{and} }j_{x}\subseteq j_{y}\text{ \emph{in} }A,
\]
\emph{where }$j_{x}$\emph{ and }$j_{y}$\emph{ are the prime ideals in }%
$A$\emph{ corresponding to }$x$\emph{ and }$y$\emph{, respectively.}

\textbf{Proof.} It is immediate from Lemmas 1.4 and 1.8.\qquad$\square$

\bigskip

Now we have got the following remark.

\bigskip

\textbf{Remark 1.10. (The Picture of a Scheme).} \emph{From the pointview of
specializations, a scheme can be regarded as a number of trees standing on the
ground such that}

$\ \left(  i\right)  $\emph{ each irreducible component is a tree;}

\emph{\ }$\left(  ii\right)  $\emph{ the initial point of an irreducible
component is the root of the corresponding tree;}

\emph{\ }$\left(  iii\right)  $\emph{ the final points of an irreducible
component are the top leaves of the corresponding tree;}

\emph{\ }$\left(  iv\right)  $\emph{ each specialization in an irreducible
component are the branches of the corresponding tree.}

\bigskip

\section{Definition for Lengths of Specializations}

\bigskip

In this section we will define the lengths of specializations, which will be
served to define the norm of a morphism of schemes in Section 3. There
exist differences between the dimensions and the lengths of subsets in a
scheme (Remarks 2.2-3).

\bigskip

Let $E$ be a topological space, and $x,y\in E$ with $x\rightarrow y.$ By a
\textbf{restrict series of specializations} from $x$ to $y$ in $E,$ denoted by
$\Gamma\left(  x,y\right)  ,$ we understand a series of specializations
\[
x=x_{0}\rightarrow x_{1}\rightarrow\cdots\rightarrow x_{n}=y
\]
in $E$ satisfying

$\ \left(  i\right)  $ $y=x$ if $x$ is a specialization of $y$ in $E;$

$\ \left(  ii\right)  $ each $x_{i}$ is not a specialization of $x_{i+1}$ for
$i=0,1,\cdots,n-1$ if $x$ is not a specialization of $y$ in $E.$

The {\textbf{{{length}}}} of the restrict series $\Gamma\left(  x,y\right)  $
is defined to be $n.$ The {{\textbf{length}}} from $x$ to $y$ (or the
{{{\textbf{length of the specialization}}}} $x\rightarrow y$), denoted by
$l\left(  x,y\right)  ,$ is defined to be the supremum among all the lengths
of restrict series of specializations from $x$ to $y$.

Set
\[
l\left(  E\right)  =\sup\{l\left(  x,y\right)  \mid x,y\in E\text{ such that
}x\rightarrow y\}.
\]
Then $l\left(  E\right)  $ is said to be the {{{\textbf{\textbf{length }of the
topological space}}}} $E.$

A restrict series $\Gamma$ of specializations in $E$ is called a
{\textbf{{{\textbf{presentation }for the length}}}} of $E$ if the length of
$\Gamma$ is equal to $l\left(  E\right)  .$ The {\textbf{{{length of a point}%
}}} $x\in E$, denoted by $l(x)$, is defined to be the length of the subspace
${Sp}(x)$ in $E$.

Obviously, $l\left(  E_{0}\right)  \leq l\left(  E\right)  $ holds for any
subspace $E_{0}$ of $E$ since a restrict series of specializations in $E_{0}$
must be in $E.$

\bigskip

\textbf{Lemma 2.1.} \emph{Let }$E$\emph{ be a topological space such that
}$\dim E<\infty.$\emph{ The following statements are true.}

$\left(  i\right)  $\emph{ Let }$\Gamma\left(  x_{0},x_{n}\right)  $\emph{ be
a presentation for the length of }$E.$\emph{ Then }$x_{0}$\emph{ is initial
and }$x_{n}$\emph{ is final in }$E$\emph{.}

$\left(  ii\right)  $\emph{ Let }$E$\emph{ be of the }$(UIP)-$\emph{property.
Then }$l\left(  E\right)  =\dim E\emph{.}$

\textbf{Proof. }$\left(  i\right)  $ It is immediate from definition.

$\left(  ii\right)  $ Hypothesize that $l\left(  E\right)  =\infty.$ That is,
for any $n\in\mathbb{N}$ there is a restrict series of specializations in $E$
\[
x_{0}\rightarrow x_{1}\rightarrow\cdots\rightarrow x_{n}.
\]
Then we have a chain of closed subsets in $E$
\[
Sp\left(  x_{0}\right)  \supsetneqq Sp\left(  x_{1}\right)  \supsetneqq
\cdots\supsetneqq Sp\left(  x_{n}\right)  .
\]

As $Sp\left(  x_{0}\right)  ,Sp\left(  x_{1}\right)  ,\cdots,Sp\left(
x_{n}\right)  $ are irreducible closed subsets, we will get $\dim E=\infty,$
which is in contradiction with the assumption that $\dim E<\infty$. Hence, we
must have $l\left(  E\right)  <\infty.$ This also proves that $l\left(
E\right)  \leq\dim E.$

Let $\dim E=n_{0}$. We have a chain of irreducible closed subsets in $E$
\[
F_{0}\supsetneqq F_{1}\supsetneqq\cdots\supsetneqq F_{n_{0}}.
\]

As $E$ is of the {{(UIP)-property}}, each subset $F_{j}$ has the unique
initial point $x_{F_{j}}$. There is a restrict series of specializations in
$E$
\[
x_{F_{0}}\rightarrow x_{F_{1}}\rightarrow\cdots\rightarrow x_{F_{n_{0}}}.
\]

Then $l\left(  E\right)  \geq\dim E$ holds. This completes the proof$.$
$\hspace{2em}\square$ {\bigskip}

{\textbf{Remark 2.2.}} \emph{Let }$E$\emph{ be a topological space of the
}$(UIP)-$\emph{property. Take a subspace }$E_{0}$\emph{ of }$E.$\emph{ In
general, it is not true that }$l\left(  E_{0}\right)  =\dim E_{0};$\emph{ but
for the whole space }$E,$\emph{ }$l\left(  E\right)  $\emph{ and }$\dim
E$\emph{ coincide with each other. That is due to the fact}

$\ \left(  i\right)  $\emph{ }$\dim E_{0}$\emph{ is determined by }$E_{0}%
$\emph{ itself;}

$\ \left(  ii\right)  $\emph{ }$l\left(  E_{0}\right)  $\emph{ is defined both
by }$E_{0}$\emph{ and by }$E$\emph{ externally.}

\bigskip

{\textbf{Remark 2.3.}}\emph{ Let }$X$\emph{ be a scheme, and }$X_{0}$\emph{ a
subscheme of }$X.$

$\left(  i\right)  $\emph{ }$l\left(  X_{0}\right)  \geq\dim X_{0}.$

$\left(  ii\right)  $\emph{ Let }$X_{0}$\emph{ be closed in }$X.$\emph{ Then
}$\dim X_{0}<\infty$\emph{ if and only if }$l\left(  X_{0}\right)  <\infty
$\emph{; moreover, }$\dim X_{0}=l\left(  X_{0}\right)  $\emph{ if }$\dim
X_{0}<\infty.$

\bigskip

\bigskip

\section{Main Results}

\bigskip

Let $x,y$ be two points in a topological space $E$. Then $x$ and $y$ are
said to be $Sp-$\textbf{connected }if either $x\rightarrow y$ in $E$ or $%
y\rightarrow x$ in $E$ holds; $x$ and $y$ are said to be $Sp-$\textbf{%
disconnected} if they are not $Sp-$connected.

\bigskip

\textbf{Definition 3.1.} \emph{Let }$f:X\rightarrow Y$\emph{\ be a
morphism
of schemes. }$\left( i\right) $\emph{\ }$f$\emph{\ is said to be \textbf{%
bounded} if there exists a constant }$\beta\in R$\emph{\ such that \ }%
\begin{equation*}
l\left( f\left( x_{1}\right) ,f\left( x_{2}\right) \right) \leq\beta l\left(
x_{1},x_{2}\right)
\end{equation*}
\emph{holds for any }$x_{1}\rightarrow x_{2}$\emph{\ in }$X$\emph{\ with }$%
l\left( x_{1},x_{2}\right) <\infty$\emph{.}

$\left( ii\right) $\emph{\ Let }$f$\emph{\ be bounded. If }$\dim X=0,$\emph{%
\ define }$\left\Vert f\right\Vert =0;$\emph{\ if }$\dim X>0,$\emph{\ define
}%
\begin{equation*}
\left\Vert f\right\Vert =\sup\{\frac{l\left( f\left( x_{1}\right) ,f\left(
x_{2}\right) \right) }{l\left( x_{1},x_{2}\right) }:x_{1}\rightarrow x_{2}%
\text{\emph{\ in }}X,\text{ }0<l\left( x_{1},x_{2}\right) <\infty\}.
\end{equation*}
\emph{Then the number }$0\leq\left\Vert f\right\Vert \leq\beta$\emph{\ is
said to be the \textbf{norm} of }$f.$

\bigskip

\textbf{Example 3.2.}

$\left( i\right) $ The $k-$rational points of a $k-$variety are morphisms of
norm zero.

$\left( ii\right) $ Let $s,t$ be varibles over a field $k.$ Suppose
\begin{equation*}
f:Spec\left( k\left[ s,t\right] \right) \rightarrow Spec\left( k\left[ t%
\right] \right)
\end{equation*}
is induced from the embedding of the $k-$algebras. Then $\left\Vert
f\right\Vert =1.$

$\left( iii\right) $ Let $t$ be a varible over $\mathbb{Q}.$ Suppose
\begin{equation*}
f:Spec\left( k\left[ t\right] \right) \rightarrow Spec\left( \mathbb{Z}\left[
t\right] \right)
\end{equation*}%
is induced from the evident embedding. Then $\left\Vert f\right\Vert =2.$

\bigskip

Here are the main results of the paper.

\bigskip

\textbf{Theorem 3.3.} \emph{Let }$f:X\rightarrow Y$\emph{\ be a morphism
of schemes. Then }$f$\emph{\ is bounded and }$0\leq \left\Vert f\right\Vert
\leq 1$ \emph{if }$f$\emph{\ satisfies }

\qquad\textbf{Condition }$\left( \ast\right) $\textbf{:}\emph{\ For any }$%
x\in X$ \emph{, either }%
\begin{equation*}
f\left( Sp\left( x\right) \right) =Sp\left( f\left( x\right) \right)
\end{equation*}
\emph{holds or }%
\begin{equation*}
f^{-1}\left( Sp\left( f\left( x\right) \right) \right) \not =\varnothing
\end{equation*}
\emph{is }$Sp-$\emph{connected.}

\bigskip

\textbf{Proof.} Without loss of generality, assume $\dim X>0$ and $\dim
Y>0.$ Let $X$ and $Y$ be both irreducible. Take any affine open subset $U$
of $X.$ We will prove
\begin{equation*}
\left\Vert f\mid_{U}\right\Vert \leq1.
\end{equation*}

Take any $x_{1}\rightarrow x_{2}$ in $U$ such
\begin{equation*}
0<\left( x_{1},x_{2}\right) <\infty\text{ and }l\left( f\left( x_{1}\right)
,f\left( x_{2}\right) \right) >0.
\end{equation*}

We will proceed in two steps.

$\left( i\right) $ Let $x_{1}\rightarrow x_{2}$ in $U$ be closest. Hypothesize
that $f\left( x_{1}\right) \rightarrow f\left( x_{2}\right) $ in $Y$ is not closest.
That is, assume
\begin{equation*}
l\left( f\left( x_{1}\right) ,f\left( x_{2}\right) \right) \geq2
\end{equation*}
in $Y.$ We have
\begin{equation*}
Sp\left( f\left( x_{1}\right) \right) \supsetneqq Sp\left( f\left(
x_{2}\right) \right) ;
\end{equation*}
then $\dim Sp\left( f\left( x_{1}\right) \right) \geq2.$

Take $y_{0}\in Y$ such that
\begin{equation*}
f\left( x_{1}\right) \rightarrow y_{0}\rightarrow f\left( x_{2}\right) \text{
in }Y
\end{equation*}
and
\begin{equation*}
f\left( x_{1}\right) \not =y_{0}\not =f\left( x_{2}\right) .
\end{equation*}
We have
\begin{equation*}
Sp\left( f\left( x_{1}\right) \right) \supsetneqq Sp\left( y_{0}\right)
\supsetneqq Sp\left( f\left( x_{2}\right) \right)
\end{equation*}
in $Y$. As
\begin{equation*}
f^{-1}\left( Sp\left( f\left( x_{1}\right) \right) \right) \supsetneqq
f^{-1}\left( Sp\left( y_{0}\right) \right) \supsetneqq f^{-1}\left( Sp\left(
f\left( x_{2}\right) \right) \right)
\end{equation*}
hold in $X$, we obtain
\begin{equation*}
\dim U\bigcap f^{-1}\left( Sp\left( f\left( x_{1}\right) \right) \right)
\geq2
\end{equation*}
from Condition $\left( \ast\right) $; then
\begin{equation*}
l\left( x_{1},x_{2}\right) \geq2,
\end{equation*}
and hence $x_{1}\rightarrow x_{2}$ in $U$ is not closest, where there will be
a contradiction. This proves $f\left( x_{1}\right) \rightarrow f\left(
x_{2}\right) $ in $Y$ is closest.

$\left( ii\right) $ Assume that $x_{1}\rightarrow x_{2}$ in $U$ is not closest.
Let
\begin{equation*}
l\left( x_{1},x_{2}\right) =n
\end{equation*}
in $U.$ Then there are the closest specializations in $U$
\begin{equation*}
z_{1}\rightarrow z_{2}\rightarrow\cdots\rightarrow z_{n+1}
\end{equation*}
where $z_{1}=x_{1}$ and $z_{n+1}=x_{2}.$

Obviously, we have either
\begin{equation*}
f\left( z_{i}\right) =f\left( z_{i+1}\right) \text{ for some }1\leq i\leq n
\end{equation*}
or

a restrict series of specializations
\begin{equation*}
f\left( z_{1}\right) \rightarrow f\left( z_{2}\right) \rightarrow
\cdots\rightarrow f\left( z_{n+1}\right)
\end{equation*}
in $Y.$

For the latter case, by $\left( i\right) $ it is seen that these specializations
are closest in $Y.$

Hence,
\begin{equation*}
l\left( f\left( x_{1}\right) ,f\left( x_{2}\right) \right) >l\left(
x_{1},x_{2}\right)
\end{equation*}
never holds in $U.$ This proves
\begin{equation*}
\left\Vert f\mid_{U}\right\Vert \leq1,
\end{equation*}
and it follows that
\begin{equation*}
\left\Vert f\right\Vert \leq1
\end{equation*}
holds in $X$ since any two $x_{1},x_{2}\in X$ with $x_{1}\rightarrow x_{2}$
are contained in an affine open subset of $X$.$\qquad\square$

\bigskip

\textbf{Remark 3.4. (The Picture of a Morphism of Schemes).
}\emph{Theorem
3.3 affords us a longitudinal classification of morphisms of schemes. Let }$%
f:X\rightarrow Y$\emph{\ be a morphism of schemes. }

$\left( i\right) $\emph{\ }$f$\emph{\ is \textbf{length-preserving} if we
have }%
\begin{equation*}
l\left( x,y\right) =l\left( f\left( x\right) ,f\left( y\right) \right)
\end{equation*}
\emph{\ for any }$x,y\in X$\emph{\ such that }$x\rightarrow y.$

$\left( ii\right) $\emph{\ }$f$\emph{\ is \textbf{asymptotic} if }$%
\left\Vert f\right\Vert =1$.

$\left( iii\right) $\emph{\ }$f$\emph{\ is \textbf{null} if }$\left\Vert
f\right\Vert =0.$

\bigskip

\textbf{Remark 3.5. (The Picture of a Morphism of Schemes). }\emph{There
exists a latitudinal classification of morphisms of schemes. That is, let }$%
f:X\rightarrow Y$\emph{\ be a morphism of schemes. }

$\left( i\right) $\emph{\ }$f$\emph{\ is \textbf{level-separated} if }$%
f\left( x\right) $\emph{\ and }$f\left( y\right) $\emph{\ are }$Sp-$\emph{%
disconnected in }$Y$\emph{\ for any }$x,y\in X$\emph{\ which are }$Sp-$\emph{%
disconnected and of the same lengths in }$X.$

$\left( ii\right) $ $f$\emph{\ is \textbf{level-reduced} if }$f\left( x\right)
$\emph{\ and }$f\left( y\right) $\emph{\ are }$Sp-$\emph{connected in
}$Y$\emph{\ for any }$x,y\in X$\emph{\ which are
}$Sp-$\emph{disconnected and of the same lengths in }$X.$

$\left( iii\right) $ $f$\emph{\ is \textbf{level-mixed} if f is neither
level-separated nor level-reduced.}

\bigskip

\textbf{Remark 3.6.} \emph{Let }$f:X\rightarrow Y$\emph{\ be a morphism
of schemes. }

$\left( i\right) $\emph{\ Let }$\dim X>0.$\emph{\ Then }$\left\Vert
f\right\Vert =1$ \emph{if }$f$\emph{\ is length-preserving.}

$\left( ii\right) $ \emph{Let }$\dim X=\dim Y<\infty .$ \emph{Then }$%
\left\Vert f\right\Vert \geq 1$ \emph{if }$f$ \emph{is surjective.}

$\left( iii\right) $\emph{\ Let }$\left\Vert f\right\Vert =1.$\emph{\ In
general, it is not true that }$f$\emph{\ is injective since there exists a
scheme }$X$\emph{\ which can not be totally ordered by specializations.}

\bigskip

\textbf{Theorem 3.7.} \emph{Let }$X$\emph{\ and }$Y$ \emph{be
irreducible
schemes, and }$f:X\rightarrow Y$ \emph{be a morphism satisfying Condition }$%
\left( \ast \right) $\emph{. Suppose }$\dim Y<\infty .$\emph{\ Then }$f$%
\emph{\ is injective if and only if }$f$\emph{\ is length-preserving and
level-separated.}

\textbf{Proof.} Prove $\implies.$ Assume that $f$ is injective. As $\dim
Y<\infty,$ we have $\dim X<\infty$ by Proposition 1.3.

Show $f$ is length-preserving. Take any restrict series of specializations in
$X$
\begin{equation*}
z_{1}\rightarrow z_{2}\rightarrow\cdots\rightarrow z_{n}.
\end{equation*}
We have
\begin{equation*}
f\left( z_{1}\right) \rightarrow f\left( z_{2}\right) \rightarrow
\cdots\rightarrow f\left( z_{n}\right)
\end{equation*}
in $Y.$ As $f$ is injective, we get $f\left( z_{i}\right) \not =f\left( z_{j}\right) $
for all $i\not =j;$ then
\begin{equation*}
l\left( z_{1},z_{n}\right) =l\left( f\left( z_{1}\right) ,f\left(
z_{n}\right) \right)
\end{equation*}
holds for any specialization $z_{1}\rightarrow z_{n}$ in $X$ which is of finite
length. As
\begin{equation*}
\dim X=l\left( X\right) <\infty,
\end{equation*}
it is seen that
\begin{equation*}
l\left( x,y\right) =l\left( f\left( x\right) ,f\left( y\right) \right)
\end{equation*}
holds for any $x\rightarrow y$ in $X.$

Show $f$ is level-separated. Take any $x_{1},x_{2}\in X$ which are $Sp-$%
disconnected and of the same lengths, that is, $x_{1}\not =x_{2}$ and $%
l\left( x_{1}\right) =l\left( x_{2}\right) ;$ then $f\left( x_{1}\right) $ and $f\left(
x_{2}\right) $ are $Sp-$disconnected; otherwise, if $f\left( x_{1}\right)
\rightarrow f\left( x_{2}\right) $ in $Y,$ we have
\begin{equation*}
Sp\left( f\left( x_{1}\right) \right) =Sp\left( f\left( x_{2}\right) \right)
\end{equation*}
since
\begin{equation*}
l\left( f\left( x_{1}\right) \right) =l\left( f\left( x_{2}\right) \right)
\leq l\left( Y\right) =\dim Y<\infty;
\end{equation*}
it follows that $f\left( x_{1}\right) =f\left( x_{2}\right) $ holds, which is in
contradiction with the assumption.

Prove $\impliedby.$ Conversely, suppose that $f$ is length-preserving and
level-separated. We have $\dim X<\infty$. In deed, if $\dim X=\infty,$ we
will obtain $\dim Y=\infty$ since $f$ is length-preserving.

Take any $x,y\in X.$ We will prove $f\left( x\right) \not =f\left( y\right) $ if
$x\not =y$.

Let $\xi$ be the generic point of $X.$ There are three cases.

\emph{Case }$\left( i\right) :$ Let $\dim X=0.$ Then $x=y$ and $f$ is
injective.

\emph{Case }$\left( ii\right) :$\emph{\ }Let $\dim X>0$ and $x=\xi.$

We have $y\neq\xi$ and $x\rightarrow y;$ then $l\left( x,y\right) >0$; as $f$
is length-preserving, it is seen that
\begin{equation*}
l\left( f\left( x\right) ,f\left( y\right) \right) =l\left( x,y\right) >0.
\end{equation*}
Hence, $f\left( x\right) \neq f\left( y\right) .$

\emph{Case }$\left( iii\right) :$ Let $\dim X>0,$ $x\neq\xi,$ and $y\neq \xi. $
As $\dim X<\infty,$ we have $l\left( z\right) \leq l\left( X\right) <\infty $ for
any $z\in X.$ There are several subcases.

If $l\left( x\right) =l\left( y\right) $ and $y\in Sp\left( x\right) $ (or $%
x\in Sp\left( y\right) ,$ respectively), we have $x=y$ and then $f\left(
x\right) =f\left( y\right) .$

If $l\left( x\right) =l\left( y\right) ,$ $y\not \in Sp\left( x\right) ,$ and $x\not
\in Sp\left( y\right) $ hold, we have $f\left( y\right) \not \in Sp\left( f\left(
x\right) \right) $ and then $f\left( x\right) \not =f\left( y\right) $ since $f$ is
level-separated.

If $l\left( x\right) >l\left( y\right) $ and $y\in Sp\left( x\right) $ hold, we have
\begin{equation*}
l\left( f\left( x\right) ,f\left( y\right) \right) =l\left( x,y\right) >0
\end{equation*}
since $f$ is length-preserving; hence, $f\left( x\right) \neq f\left( y\right) .$

Now suppose $l\left( x\right) >l\left( y\right) $ and $y\notin Sp\left(
x\right) $ without loss of generality. It is seen that $x,y$ are $Sp-$%
disconnected. Taking a presentation for the length $l\left( x\right) $, we
have $x_{0}\in Sp\left( x\right) $ such that $\left. l\left( x_{0}\right) =l(y\right)
<\infty$. Then $l\left( x,x_{0}\right) >0$ and $x_{0}\neq y.$ As $f$ is
level-separated, we have $f\left( x_{0}\right) \neq f\left( y\right) . $

As $l\left( x_{0}\right) =l\left( y\right) ,$ we have $l\left( \xi ,x_{0}\right)
=l\left( \xi,y\right) $ by taking presentations for the lengths; as
\begin{equation*}
l\left( \xi,x_{0}\right) =l\left( \xi,x\right) +l\left( x,x_{0}\right)
\end{equation*}
holds, we have
\begin{equation*}
l\left( \xi,y\right) =l\left( \xi,x\right) +l\left( x,x_{0}\right) .
\end{equation*}
Hence,
\begin{equation*}
l\left( f\left( x\right) ,f\left( x_{0}\right) \right) =l\left(
x,x_{0}\right) >0
\end{equation*}
and
\begin{equation*}
l\left( f\left( \xi\right) ,f\left( y\right) \right) =l\left( f\left(
\xi\right) ,f\left( x_{0}\right) \right) =l\left( f\left( \xi\right)
,f\left( x\right) \right) +l\left( f\left( x\right) ,f\left( x_{0}\right)
\right) <\infty.
\end{equation*}

We must have $f\left( x\right) \neq f\left( y\right) ;$ otherwise, if $%
f\left( x\right) =f\left( y\right) ,$ there is
\begin{equation*}
l\left( f\left( \xi\right) ,f\left( x\right) \right) =l\left( f\left(
\xi\right) ,f\left( x\right) \right) +l\left( f\left( x\right) ,f\left(
x_{0}\right) \right) ;
\end{equation*}
then $l\left( f\left( x\right) ,f\left( x_{0}\right) \right) =0,$ where there will be
a contradiction. {\hspace{2em}} $\square$

\bigskip

\textbf{Corollary 3.8.} \emph{Let }$f:X\rightarrow Y$\emph{\ be a morphism
of schemes. Then }$\left\Vert f\right\Vert =1$\emph{\ if }$f$\emph{\ is
injective and satisfies Condition }$\left( \ast\right) $\emph{.}

\textbf{Proof.} As every ideal is contained in a maximal ideal in a
commutative ring, we can take an irreducible open subspace $X_{0}$ of $X$
such that $l\left( X_{0}\right) <\infty.$ Then we have $\left\Vert
f\mid_{X_{0}}\right\Vert =1;$ as $\left\Vert f\right\Vert \leq1$ by Theorem
3.3, we get $\left\Vert f\right\Vert =1.$\qquad$\square$

\bigskip

\bigskip

\section{An Application of Specializations}

\bigskip

\textbf{Definition 4.1. }$\left( i\right) $ \emph{Let }$R,S$\emph{\ be
commutative rings with }$1$\emph{. A homomorphism }$\tau:R\rightarrow S$%
\emph{\ is said to be of }$J-$\textbf{\ type}\emph{\ if }$\tau^{-1}\left(
\tau\left( I\right) S\right) =I$\emph{\ holds for every prime ideal }$I$%
\emph{\ in }$R.$

$\left( ii\right) $ \emph{Let }$X,Y$\emph{\ be schemes. A morphism }$%
f:X\rightarrow Y$\emph{\ is said to be of }\textbf{finite }$J-$\textbf{\ type%
}\emph{\ if }$f$\emph{\ is of finite type and the induced homomorphism }%
\begin{equation*}
f^{\#}\mid_{V}:O_{Y}\left( V\right) \rightarrow f_{\ast}O_{X}\left( U\right)
\end{equation*}
\emph{\ is of }$J-$\emph{type for any affine open sets }$V$\emph{\ of }$Y$%
\emph{\ and }$U$\emph{\ of }$f^{-1}\left( V\right) .$

\bigskip

\textbf{Proposition 4.2.} \emph{Let }$X$\emph{\ and }$Y$ \emph{be
irreducible schemes, and }$f:X\rightarrow Y$ \emph{be a morphism. Then
we have }$\dim X=\dim Y$\emph{\ if }$f$\emph{\ is length-preserving and
of finite }$J-$\emph{type.}

\textbf{Proof.} Let $f$ be length-preserving and of finite $J-$type. It follows
that
\begin{equation*}
l\left( X\right) =l\left( f\left( X\right) \right) \leq l\left( Y\right)
\end{equation*}
hold. Then we have $\dim X\leq\dim Y$ since
\begin{equation*}
\dim X=l\left( X\right) \text{ and }l\left( Y\right) =\dim Y.
\end{equation*}

Take any $x\in X$ and $y=f\left( x\right) \in Y.$ As $f$ is of finite $J-$%
type, there are affine open subsets $V$ of $Y$ and $U$ of $f^{-1}\left(
V\right) $ such that
\begin{equation*}
f^{\#}\mid_{V}:O_{Y}\left( V\right) \rightarrow f_{\ast}O_{X}\left( U\right)
\end{equation*}
is a homomorphism of $J-$type, where $x\in U$ and $y\in V.$

Set $V=Spec\left( R\right) $ and $U=Spec\left( S\right) .$ As $X$ and $Y$
are irreducible, we have
\begin{equation*}
\dim U=\dim X\text{ and }\dim V=\dim Y.
\end{equation*}

Take any restrict series of specializations
\begin{equation*}
y_{0}\rightarrow y_{1}\rightarrow\cdots\rightarrow y_{n}
\end{equation*}
in $V.$ Then we obtain a chain of prime ideals
\begin{equation*}
j_{y_{0}}\subsetneqq j_{y_{1}}\subsetneqq\cdots\subsetneqq j_{y_{n}}
\end{equation*}
in $R,$ where each $j_{y_{i}}$ is the prime ideal in $R$ corresponding to $%
y_{i}$ in $V.$ By Corollary 2.3$^{\left[ 5\right] }$ there are prime ideals
\begin{equation*}
I_{0}\subseteq I_{1}\subseteq\cdots\subseteq I_{n}
\end{equation*}
in $S$ such that $f^{\#-1}\left( I_{i}\right) =j_{y_{i}}.$

Hence, we obtain a restrict series of specializations
\begin{equation*}
x_{0}\rightarrow x_{1}\rightarrow\cdots\rightarrow x_{n}
\end{equation*}
in $U$ such that $f\left( x_{i}\right) =y_{i}$ and $j_{x_{i}}=I_{i}.$ This proves
$l\left( U\right) \geq l\left( V\right) .$

As
\begin{equation*}
\dim X=l\left( X\right) \geq l\left( U\right)
\end{equation*}
and
\begin{equation*}
l\left( V\right) \geq\dim V=\dim Y;
\end{equation*}
we have
\begin{equation*}
\dim X\geq\dim Y.
\end{equation*}

This completes the proof.$\qquad\square$


\newpage

\begin{thebibliography} {9}

\bibitem{1} Grothendieck, A. $\acute{E}$l$\acute{e}$ments de
G$\acute{e}$oem$\acute{e}$trie Alg$\acute{e}$brique (EGA1). IHES, 1960.

\bibitem{2} Grothendieck, A. $\acute{E}$l$\acute{e}$ments de
G$\acute{e}$om$\acute{e}$trie Alg$\acute{e}$brique (EGA2). IHES, 1961.

\bibitem{3} Grothendieck, A. $\acute{E}$l$\acute{e}$ments de
G$\acute{e}$om$\acute{e}$trie Alg$\acute{e}$brique (EGA1/4). IHES, 1964.

\bibitem{4} Grothendieck, A. $\acute{E}$l$\acute{e}$ments de
G$\acute{e}$om$\acute{e}$trie Alg$\acute{e}$brique (EGA2/4). IHES, 1965.

\bibitem{S} Sharma, P K. A note on lifting of chains of prime ideals.
Journal of Pure and Applied Algebra, 192(2004),287-291.

\bibitem{W} Weil, A. Foundations of Algebraic Geometry. Amer Math
Society, 1946.
\end{thebibliography}
\end{document}